\documentclass[11pt,twoside,reqno]{amsart}
\usepackage{amscd,amssymb}
\usepackage[arrow,matrix]{xy}

\calclayout

\newtheorem{theorem}{\bf Theorem}
\newtheorem{lemma}[theorem]{\bf Lemma}

\newtheorem{cor}[theorem]{\bf Corollary}

\newtheorem{rk}{\bf Remark}
\newtheorem{ex}{\bf Example}

\newcommand{\R}{\mathbb{R}}
\newcommand{\C}{\mathbb{C}}

\newcommand{\OO}{\mathcal{O}}

\newcommand{\ir}{\mathcal{R}}

\newcommand{\jj}{\mathcal{J}}

\begin{document}

\vspace{1.3cm}

\title
{Weighted Homogeneous Polynomials with Isomorphic Milnor Algebras}

\author{Imran Ahmed$^1$}
\thanks{ {\enskip
  \enskip $^{1}$Department of Mathematics, COMSATS Institute of Information Technology, M.A. Jinnah Campus, Defence Road, off Raiwind Road Lahore, PAKISTAN and Departamento de Matem\'{a}ticas, Instituto de Ci\^{e}ncias Matem\'{a}ticas e
de Computa\c{c}\~{a}o, Universidade de S\~{a}o Paulo, Avenida Trabalhador
S\~{a}ocarlense 400, S\~{a}o
Carlos-S.P., Brazil. Email: drimranahmed@ciitlahore.edu.pk}}
\begin{abstract}
We recall first some basic facts on weighted homogeneous functions and filtrations in the ring $A$ of formal power series. We introduce next their analogues for weighted homogeneous diffeomorphisms and vector fields. We show that the Milnor algebra is a complete invariant for the classification of weighted homogeneous polynomials with respect to right-equivalence, i.e. change of coordinates in the source and target by diffeomorphism. \vskip 0.4 true cm
 \noindent
 {\it AMS SUBJECT} : Primary 14L30, 14J17, 16W22.\\
\end{abstract}
\maketitle


\pagestyle{myheadings} \markboth{\centerline {\scriptsize Imran Ahmed}}
         {\centerline {\scriptsize  Weighted Homogeneous Polynomials with Isomorphic Milnor Algebras}}


\bigskip
\bigskip
\medskip

\section{Introduction}
Let $f\in\C[x_1,\ldots,x_n]$ be a weighted homogeneous polynomial of
degree $d$ w.r.t weights $(w_1,\ldots,w_n)$ and $f_i=\frac{\partial
f}{\partial x_i},\,\,i=1,\ldots,n$ its partial derivatives. The
Milnor algebra of $f$ is defined by
$$M(f)=\frac{\C[x_1,\ldots,x_n]}{\jj_f},$$ where $\jj_f=\langle f_1,\ldots,f_n\rangle$ is the Jacobian ideal.

We say that two weighted homogeneous polynomials $f,g:\C^n\to \C$ are $\ir$-equivalent if there exists a diffeomorphism $\psi:\C^n\to\C^n$ such that $f\circ\psi=g$.

We recall first some basic facts on weighted homogeneous functions and filtrations in the ring $A$ of formal power series. We introduce next their analogues for weighted homogeneous diffeomorphisms and vector fields.

In Theorem \ref{th8.1} we show that two weighted homogeneous polynomials $f$ and $g$ having isomorphic Milnor algebras are right-equivalent. The Example of Gaffney and Hauser, in \cite{GH}, suggests us that we
can not extend this result for arbitrary analytic germs.

\section{Preliminary Results}
We recall first some basic facts on weighted homogeneous functions and filtrations in the ring $A$ of formal power series. We introduce next their analogues for weighted homogeneous diffeomorphisms and vector fields. For a more complete introduction see \cite{Ar}, Chap. 1, \S 3.

A holomorphic function $f:(\C^n,0)\to (\C,0)$ (defined on the complex space $\C^n$) is a weighted homogeneous function of degree $d$ with weights $w_1,\ldots, w_n$ if
$$f(\lambda^{w_1}x_1,\ldots,\lambda^{w_n}x_n)=\lambda^d f(x_1,\ldots,x_n)\,\forall\, \lambda >0.$$

In terms of the Taylor series $\sum f_{\underline{k}}x^{\underline{k}}$ of $f$, the weighted homogeneity
condition means that the exponents of the nonzero terms of the series lie in the hyperplane
$$L=\{\underline{k}:w_1k_1+\ldots+w_n k_n=d\}.$$

Any weighted homogeneous function $f$ of degree $d$ satisfies Euler's identity
\begin{equation}\label{euler}
\sum_{i=1}^n w_i x_i \frac{\partial f}{\partial x_i}=d.f
\end{equation}
It implies that a weighted homogeneous function $f$ belongs to its Jacobean ideal $\jj_f$. The necessary and sufficient conditions of a function-germ $f:(\C^n,0)\to (\C,0)$ to be equivalent to a weighted homogeneous function-germ are $f\in\jj_f$, which is the well known result of Saito \cite{Sa}.

Consider $\C^n$ with a fixed coordinate system $x_1,\ldots,x_n$. The
algebra of formal power series in the coordinates will be denoted by
$A=\C[[x_1,\ldots,x_n]]$. We assume that a weighted homogeneity type
$\underline{w}=(w_1,\ldots,w_n)$ is fixed. With each such
$\underline{w}$ there is associated a filtration of the ring $A$,
defined as follows.

The monomial ${\bf x}^{\underline{k}}$ is said to have degree $d$ if
$<\underline{w},\underline{k}>=w_1k_1+\ldots+w_n k_n=d$.

The order $d$ of a series (resp. polynomial) is the smallest of the
degrees of the monomials that appear in that series (resp.
polynomial).

The series of order larger than or equal to $d$ form a subspace
$A_d\subset A$. The order of a product is equal to the sum of the
orders of the factors. Consequently, $A_d$ is an ideal in the ring
$A$. The family of ideals $A_d$ constitutes a decreasing filtration
of $A$: $A_{\acute{d}}\subset A_d$ whenever $\acute{d}>d$. We let
$A_{d+}$ denote the ideal in $A$ formed by the series of order
higher than $d$.

The quotient algebra $A/A_{d+}$ is called the algebra of
$d$-weighted jets, and its elements are called $d$-weighted jets.

\medskip

Several Lie groups and algebras are associated with the filtration
defined in the ring $A$ of power series by the type of
weighted homogeneity $\underline{w}$. In the case of ordinary
homogeneity these are the general linear group, the group of
$k$-jets of diffeomorphisms, its subgroup of $k$-jets with
$(k-1)$-jet equal to the identity, and their quotient groups. Their
analogues for the case of a weighted homogeneous filtration are defined
as follows.

A formal diffeomorphism $g:(\C^n,0)\to (\C^n,0)$ is a set of $n$
power series $g_i\in A$ without constant terms for which the map
$g^{\ast}:A\to A$ given by the rule $g^{\ast}f=f\circ g$ is an
algebra isomorphism.

The diffeomorphism $g$ is said to have order $d$ if for every $s$
$$(g^{\ast}-1)A_s\subset A_{s+d}.$$

The set of all diffeomorphisms of order $d\geq0$ is a group $G_d$.
The family of groups $G_d$ yields a decreasing filtration of the
group $G$ of formal diffeomorphisms; indeed, for $\acute{d}>d\geq0$,
$G_{\acute{d}}\subset G_d$ and is a normal subgroup in $G_d$.

The group $G_0$ plays the role in the weighted homogeneous case that the
full group of formal diffeomorphisms plays in the homogeneous case.
We should emphasize that in the weighted homogeneous case $G_0\neq G$,
since certain diffeomorphisms have negative orders and do not belong
to $G_0$.

The group of $d$-weighted jets of type $\underline{w}$ is the quotient
group of the group of diffeomorphisms $G_0$ by the subgroup $G_{d+}$
of diffeomorphisms of order higher than $d$: $J_d=G_0/G_{d+}$.

Note that in the ordinary homogeneous case our numbering differs from the
standard one by $1$: for us $J_0$ is the group of $1$-jets and so
on.

$J_d$ acts as a group of linear transformations on the space $A/A_{d+}$ of $d$-weighted jets of functions. A special importance is attached to the group $J_0$, which is the weighted homogeneous generalization of the general linear group.

A diffeomorphism $g\in G_0$ is said to be weighted homogeneous of type
$\underline{w}$ if each of the spaces of weighted homogeneous functions
of degree $d$ (and type \underline{w}) is invariant under the action
of $g^{\ast}$.

The set of all weighted homogeneous diffeomorphisms is a subgroup of
$G_0$. This subgroup is canonically isomorphic to $J_0$, the
isomorphism being provided by the restriction of the canonical
projection $G_0\to J_0$.

\medskip

The infinitesimal analogues of the concepts introduced above look as
follows.

A formal vector field $v=\sum v_i\partial_i$, where
$\partial_i=\partial/\partial x_i$, is said to have order $d$ if
differentiation in the direction of $v$ raises the degree of any
function by at least $d$: $L_v A_s\subset A_{s+d}$.

We let $\mathfrak{g}_d$ denote the set of all vector fields of order $d$. The
filtration arising in this way in the Lie algebra $\mathfrak{g}$ of vector
fields (i.e., of derivations of the algebra $A$) is compatible with
the filtrations in $A$ and in the group of diffeomorphisms $G$:\\
1. $f\in A_d, v\in \mathfrak{g}_s\Rightarrow fv\in \mathfrak{g}_{d+s}, L_vf\in A_{d+s}$\\
2. The module $\mathfrak{g}_d$, $d\geq 0$, is a Lie algebra w.r.t. the Poisson
bracket of vector fields.\\
3. The Lie algebra $\mathfrak{g}_d$ is an ideal in the Lie algebra $\mathfrak{g}_0$.\\
4. The Lie algebra $\mathfrak{j}_d$ of the Lie group $J_d$ of $d$-weighted jets of
diffeomorphisms is equal to the quotient algebra $\mathfrak{g}_0/\mathfrak{g}_{d+}$.\\
5. The weighted homogeneous vector fields of degree $0$ form a finite-dimensional Lie subalgebra of the Lie algebra $\mathfrak{g}_0$; this subalgebra
is canonically isomorphic to the Lie algebra $\mathfrak{j}_0$ of the group of
$0$-jets of diffeomorphisms.

The support of a weighted homogeneous function of degree $d$ and type $\underline{w}$ is the set of all points $\underline{k}$ with nonnegative integer coordinates on the diagonal $$L=\{\underline{k}:\langle\underline{k},\underline{w}\rangle=d\}.$$

Weighted homogeneous functions can be regarded as functions given on their supports: $\sum f_{\underline{k}}x^{\underline{k}}$ assumes at $\underline{k}$ the value $f_{\underline{k}}$. The set of all such functions is a linear space $\C^r$, where $r$ is the number of points in the support. Both the group of weighted homogeneous diffeomorphisms (of type \underline{w}) and its Lie algebra $\mathfrak{a}$ act on this space.

The Lie algebra $\mathfrak{a}$ of a weighted homogeneous vector field of degree $0$
is spanned, as a $\C$-linear space, by all monomial fields
$x^{\underline{P}}\partial_i$ for which
$<\underline{P},\underline{w}>=w_i$. For example, the $n$ fields
$x_i\partial_i$ belong to $\mathfrak{a}$ for any $\underline{w}$.

\begin{ex}
Consider the weighted homogeneous polynomial $f=x^2y+z^2$ of degree
$d=6$ w.r.t. weights $(2,2,3)$. Note that the Lie algebra of
weighted homogeneous vector fields of degree $0$ is spanned by
$$\mathfrak{a}=\langle x^{\underline{P}}\partial_i:<\underline{P},\underline{w}>=w_i, i=1,2,3 \rangle
=\langle x\frac{\partial}{\partial x},x\frac{\partial}{\partial
y},y\frac{\partial}{\partial x}, y\frac{\partial}{\partial
y},z\frac{\partial}{\partial z}\rangle$$
\end{ex}

\section{Main Results}

We recall first Mather's lemma providing effective necessary  and sufficient conditions for a connected submanifold (in our case the path $P$) to be contained in an orbit.

\begin{lemma} (\cite{Ma})\label{lem4.3.1}
Let $m:G\times M\to M$ be a smooth action and $P\subset M$ a
connected smooth submanifold. Then $P$ is contained in a single
$G$-orbit if and only if the following conditions are
fulfilled:\\
(a) $T_x(G.x)\supset T_xP$, for any $x \in P$.\\
(b) $\dim T_x(G.x)$ is constant for $x\in P$.
\end{lemma}

For arbitrary (i.e. not necessary with isolated singularities) weighted homogeneous polynomials we establish the following result.

\begin{theorem}\label{th8.1}
Let $f,g$ be two weighted homogeneous polynomials of degree $d$ w.r.t. weights
$(w_1,\ldots,w_n)$ such that $\jj_f=\jj_g$. Then
$f\stackrel{\ir}{\sim}g$, where $\stackrel{\ir}{\sim}$ denotes the
right equivalence.
\end{theorem}

\begin{proof}
Let $H_{\underline{w}}^d(n,1;\C)$ be a space of weighted homogeneous polynomials from $\C^n$ to $\C$ of degree $d$ w.r.t weights $(w_1,\ldots,w_n)$. Let $f,g\in H_{\underline{w}}^d(n,1;\C)$ such that $\jj_f=\jj_g$. Set $f_t=(1-t)f+tg\in H_{\underline{w}}^d(n,1;\C)$. Consider the $\ir$-equivalence action on $H_{\underline{w}}^d(n,1;\C)$ under the
group of $1$-jets $J_0$, we have
\begin{equation}\label{eq8.1}
T_{f_t}(J_0.f_t)=\C\langle x^{\underline{P}}\frac{\partial
f_t}{\partial
x_i};\,\,i=1,\ldots,n\,\mbox{and}\,<\underline{P},\underline{w}>=w_i\rangle
\end{equation}
Note that $T_{f_t}(J_0.f_t)\subset\jj_{f_t}\cap H_{\underline{w}}^d$.
But $\jj_{f_t}\cap H_{\underline{w}}^d\subset \jj_f\cap
H_{\underline{w}}^d$ since
\[\frac{\partial f_t}{\partial x_i}=(1-t)\frac{\partial f}{\partial
x_i}+t\frac{\partial g}{\partial
x_i}\in(1-t)\jj_f+t\jj_g=\jj_f\,\,\,\,\,\,\,\,(\mbox{because }\jj_f=\jj_g)\] Therefore, we
have the inclusion of finite dimensional $\C$-vector spaces
\begin{equation}\label{eq8.2}
T_{f_t}(J_0.f_t)=\C\langle x^{\underline{P}}\frac{\partial
f_t}{\partial
x_i};\,\,i=1,\ldots,n\,\mbox{and}\,<\underline{P},\underline{w}>=w_i\rangle\subset
\jj_f\cap H_{\underline{w}}^d
\end{equation}
with equality for $t=0$ and $t=1$.\\Let's  show that we have
equality for all $t\in[0,1]$ except finitely many values.\\
Take $\dim(\jj_f\cap H_{\underline{w}}^d)=m$ (say). Let's fix $\{e_1,\ldots,e_m\}$ a basis of
$\jj_f\cap H_{\underline{w}}^d$. Consider the
$m$ polynomials corresponding to the generators of the space
\eqref{eq8.1}:
\[\alpha_i(t)=x^{\underline{P}}\frac{\partial f_t}{\partial
x_i}=x^{\underline{P}}[(1-t)\frac{\partial f}{\partial
x_i}+t\frac{\partial g}{\partial
x_i}],\,\,\mbox{where}\,<\underline{P},\underline{w}>=w_i\,\mbox{and}\,\underline{P}=(P_1,\ldots,P_n)\]

We can express each $\alpha_i(t)$, $i=1,\ldots,m$ in terms of
above mentioned fixed basis as
\begin{equation}\label{eq8.3}
\alpha_i(t)=\phi_{i1}(t)e_1+\ldots+\phi_{im}(t)e_m,\,\,\forall\,\,
i=1,\ldots,m
\end{equation}
where each $\phi_{ij}(t)$ is linear in $t$. Consider the matrix of
transformation corresponding to the eqs. \eqref{eq8.3}
\[
(\phi_{ij}(t))_{m\times m}= \left(
  \begin{array}{cccc}
    \phi_{11}(t) & \phi_{12}(t)& \ldots & \phi_{1m}(t) \\
    \vdots & \vdots & \ddots & \vdots \\
    \phi_{m1}(t) & \phi_{m2}(t) & \ldots & \phi_{mm}(t) \\
  \end{array}
\right)
\]
having rank at most $m$. Note that the equality $$\C\langle
x^{\underline{P}}\frac{\partial f_t}{\partial
x_i};\,\,i=1,\ldots,n\,\mbox{and}\,<\underline{P},\underline{w}>=w_i\rangle=\jj_f\cap
H_{\underline{w}}^d$$ holds for those values of $t$ in $\C$ for
which the rank of above matrix is precisely $m$. We have the $m\times
m$-matrix whose determinant is a polynomial of degree $m$ in $t$
and by the fundamental theorem of algebra it has at most $m$ roots
in $\C$ for which rank of the matrix of transformation will be less
than $m$. Therefore, the above-mentioned equality does not hold for
at most finitely many values, say $t_1,\ldots,t_q$ where
$1\leq q\leq m$.\\
It follows that the dimension of the space \eqref{eq8.1} is constant
for all $t\in \C$ except finitely many values $\{t_1,\ldots,t_q\}$.\\
For an arbitrary smooth path
\[\alpha:\C\longrightarrow\C\backslash \{t_1,\ldots,t_q\}\]
with $\alpha(0)=0$ and $\alpha(1)=1$, we have the connected smooth
submanifold
\[P=\{f_t=(1-\alpha(t))f(x)+\alpha(t)g(x):\,t\in\C\}\]
of $H_{\underline{w}}^d$. By the above, it follows $\dim
T_{f_t}(J_0.f_t)$ is
constant for $f_t\in P$.\\
Now, to apply Mather's lemma, we need to show that the tangent space
to the submanifold $P$ is contained in that to the orbit
$J_0.f_t$ for any $f_t\in P$. One clearly has
\[T_{f_t}P=\{\dot{f_t}=-\dot{\alpha}(t)f(x)+\dot{\alpha}(t)g(x):\,\forall\,t\in\C\}\]
Therefore, by Euler formula \ref{euler}, we
have\[T_{f_t}P\subset T_{f_t}(J_0.f_t)\] By Mather's lemma the
submanifold $P$ is contained in a single orbit. Hence the result.
\end{proof}

\begin{cor}\label{cor8.2}
Let $f,g$ be two weighted homogeneous polynomials of degree $d$ w.r.t. weights
$(w_1,\ldots,w_n)$. If $M(f)\simeq M(g)$ (isomorphism of graded $\C$-algebra) then $f\stackrel{\ir}{\sim}g$.
\end{cor}

\begin{proof}
We show firstly that an isomorphism of graded $\C$-algebras
\[\varphi:M(g)=\frac{\C[x_1,\ldots,x_n]}{\jj_g}\stackrel{\simeq}{\longrightarrow}
M(f)=\frac{\C[x_1,\ldots,x_n]}{\jj_f}\] is induced by an isomorphism
$u:\C^n\longrightarrow\C^n$ such that $u^*(\jj_g)=\jj_f$.\\Consider
the following commutative diagram.
$$\xymatrix{
0 \ar[d] & 0 \ar[d] \\
\jj_g  \ar@{.>}[r]^{u^*} \ar[d]^i & \jj_f \ar[d]^j \\
\C[x_1,\ldots,x_n]  \ar@{.>}[r]^{u^*} \ar[d]^p & \C[x_1,\ldots,x_n]
\ar[d]^q
\\
M(g) \ar[r]^{\varphi}_{\simeq} \ar[d] & M(f) \ar[d] \\
0 & 0 }
$$

Define the morphism $u^*:\C[x_1,\ldots,x_n]\rightarrow
\C[x_1,\ldots,x_n]$ by
\begin{equation}\label{eq8.4}
u^*(x_i)=L_i(x_1,\ldots,x_n)=\sum_{j=1}^n a_{ij}x_j^{\alpha_j}+\sum
a_{ik_1\ldots k_n}x_{k_1}^{\beta_1}\ldots
x_{k_n}^{\beta_n}\,;\,i=1,\ldots,n
\end{equation}
where
$k_l\in\{1,\ldots,n\}\,\mbox{and}\,w_{k_1}\beta_1+\ldots+w_{k_n}\beta_n=deg_{\underline{w}}(x_i)=w_j\alpha_j$,
which is well defined by commutativity of diagram below.
$$\xymatrix{
x_i  \ar@{|.>}[r]^{u^*} \ar@{|->}[d]^p &  L_i\ar@{|->}[d]^q \\
\widehat{x_i} \ar@{|->}[r]^{\varphi}_{\simeq} & \widehat{L_i}}$$

Note that the isomorphism $\varphi$ is a degree preserving map and
is also given by the same morphism $u^*$. Therefore, $u^*$ is an
isomorphism.

Now, we show that $u^*(\jj_g)=\jj_f$. For every $G\in\jj_g$, we
have $u^*(G)\in\jj_f$ by commutative diagram below.
$$\xymatrix{
G  \ar@{|->}[r]^{u^*} \ar@{|->}[d]^p &  F=u^*(G) \ar@{|->}[d]^q \\
\widehat{0} \ar@{|->}[r]^{\varphi} & \widehat{F}=\widehat{0}}
$$
It implies that $u^*(\jj_g)\subset\jj_f$. As $u^*$ is an
isomorphism, therefore it is invertible and by repeating the above
argument for its inverse, we have $u^*(\jj_g)\supset\jj_f$.

Thus, $u^*$ is an isomorphism with $u^*(\jj_g)=\jj_f$. By eq. \eqref{eq8.4}, the map $u:\C^n\to\C^n$ can be defined by
$$u(z_1,\ldots,z_n)=(L_1(z_1,\ldots,z_n),\ldots,L_n(z_1,\ldots,z_n))$$
where $L_i(z_1,\ldots,z_n)=\sum_{j=1}^n a_{ij}x_j^{\alpha_j}+\sum
a_{ik_1\ldots k_n}x_{k_1}^{\beta_1}\ldots
x_{k_n}^{\beta_n}\,;\,i=1,\ldots,n,\\
k_l\in\{1,\ldots,n\}\,\mbox{and}\,w_{k_1}\beta_1+\ldots+w_{k_n}\beta_n=deg_{\underline{w}}(x_i)=w_j\alpha_j$. Note that $u$ is an isomorphism by Prop. 3.16 see \cite{D1}, p.23.

In this way, we have shown that the isomorphism $\varphi$ is induced
by the isomorphism $u:\C^n\to\C^n$ such that
$u^*(\jj_g)=\jj_f$.

Consider $u^*(\jj_g)=<g_1\circ u,\ldots,g_n\circ u>=\jj_{g\circ u}$,
where $g_j=\frac{\partial g}{\partial x_j}$. Therefore, $\jj_{g\circ
u}=\jj_f\Rightarrow g\circ u\stackrel{\ir}{\sim}f$, by Theorem
\ref{th8.1}. But $g\circ u\stackrel{\ir}{\sim}g$. It follows that
$g\stackrel{\ir}{\sim}f$.
\end{proof}

\begin{rk}\label{rk8.3}
The converse implication, namely
\[f\stackrel{\ir}{\sim}g\Rightarrow M(f)\simeq M(g)\]
always holds(even for analytic germs $f,\,g$ defining IHS), see
\cite{D1}, p90.
\end{rk}

The following Example of Gaffney and Hauser \cite{GH}, suggests
us that we can not extend the Theorem \ref{th8.1} for arbitrary
analytic germs.

\begin{ex}
Let $h:(\C^n,0)\to (\C,0)$ be any function satisfying $h\notin
\jj_h\subseteq\OO_n$ i.e. $h\notin H_{\underline{w}}^d(n,1;\C)$.
Define a family $f_t:(\C^n\times\C^n\times\C,0)\to(\C,0)$ by
$f_t(x,y,z)=h(x)+(1+z+t)h(y)$, and let
$(X_t,0)\subseteq(\C^{2n+1},0)$ be the hypersurface defined by
$f_t$. Note that
$$\jj_{f_t}=\langle\frac{\partial h}{\partial x_i}(x),\frac{\partial
h}{\partial y_j}(y),h(y)\rangle,\,\,t\in\C.$$ On the other hand, the
family $\{(X_t,0)\}_{t\in\C}$ is not trivial i.e.
$(X_t,0)\ncong(X_0,0)$: For, if $\{f_t\}_{t\in\C}$ were trivial, we
would have by Proposition 2, $\S 1$, \cite{GH}
$$\frac{\partial
f_t}{\partial t}=h(y)\in
(f_t)+m_{2n+1}\jj_{f_t}=(f_t)+m_{2n+1}\jj_{h(x)}+m_{2n+1}\jj_{h(y)}+m_{2n+1}(h(y))$$
Solving for $h(y)$ implies either $h(y)\in \jj_{h(y)}$ or $h(x)\in
\jj_{h(x)}$ contradicting the assumption on $h$.\\
It follows that $f_t$ is not $\ir$-equivalent to $f_0$.
\end{ex}

\medskip
\bigskip

{\bf Acknowledgements.}  I would like to thank Maria Aparecida Soares Ruas for helpful conversations in developing this paper. This work was partially supported by FAPESP, grant \# 2010/01895-3, and CNPq-TWAS, grant \# FR 3240188100.

\end{document}